\documentclass[a4paper]{article}
\usepackage{amsmath,amssymb,amsfonts}

\newtheorem{theorem}{Theorem}[section]

\title{On Goldbach's Conjecture }

\author{Jailton C. Ferreira}

\date{ }
\begin{document}
\maketitle
\pagenumbering{arabic}

\begin{abstract}
It is shown that if every odd integer $n > 5$ is the sum of three
primes, then every even integer $n > 2$ is the sum of two primes.
A conditional proof of Goldbach's conjecture, based on Cram\'er's
conjecture, is presented. Theoretical and experimental results
available on Goldbach's conjecture allow that a less restrictive
conjecture than Cram\'er's conjecture be used in the conditional
proof. A basic result of the Maier's paper on Cram\'er's model is
criticized.
\end{abstract}

\section{Introduction} \label{sec-1}

\hspace{22pt} In 1742 Goldbach wrote a letter to Euler
conjecturing that every integer greater than 2 is the sum of three
prime numbers. Euler replied that this conjecture breaks up into
two: every even integer is the sum of two primes; every odd
integer is the sum of three primes. The conjecture \textit{``every
even integer $n > 2$ is the sum of two primes"} is now known as
Goldbach's conjecture and the conjecture \textit{``every odd
integer $n > 5$ is the sum of three primes"} is known as ``the
weaker", ``the odd" or ``the ternary" Goldbach's conjecture.

\hspace{22pt} The ternary Goldbach's conjecture, abbreviated here
as ``ternary GC", is considered the easiest of the two cases. In
1937 Vinogradov ~\cite{Vinogradov} proved that the ternary GC is
true for sufficiently large odd number. In 1956 Borodzkin
~\cite{Borodzkin} showed that odd numbers greater than
$3^{3^{15}}$ are sufficiently large in Vinogradov's proof. This
bound was reduced to $e^{e^{11503}} (\approx 10^{43000})$ by Chen
and Wang ~\cite{Chen&Wang1} in 1989 and to $e^{e^{9715}} (\approx
10^{7194})$ in 1996 ~\cite{Chen&Wang2}.

\hspace{22pt} The Goldbach's conjecture is known to be true up to
$10^{16}$. Deshouillers, te Riele and Saouter
~\cite{Deshouillers&Riele&Saouter1} have checked it up to
$10^{14}$, Richstein ~\cite{Richstein} up to $4 \times 10^{14}$
and Silva ~\cite{Silva} up to $1 \times 10^{16}$.

\hspace{22pt} Deshouillers \textit{et alli}
~\cite{Deshouillers&alli1} outlined a proof in which if the
Generalized Riemann Hypothesis holds, then the ternary GC is true.
As far as we are concerned, there is not an analogous conditional
proof for Goldbach's conjecture.

\section{The same truth value}
\label{sec-2}

\hspace{22pt} Let us notice that if Goldbach's conjecture is true
then the ternary GC is true. In the case Goldbach's conjecture is
true if the ternary GC is true, the conditional proof of
Deshouillers \textit{et alli} ~\cite{Deshouillers&alli1} can be
used to conditionally prove Goldbach's conjecture.

\begin{theorem} \label{teorema-1}
If every odd integer $n > 5$ is the sum of three primes, then
every even integer $n > 2$ is the sum of two primes.
\end{theorem}
\textit{Proof:}

\hspace{22pt}Let us assume that exists an even integer $m$ greater
than 2 that can not be expressed as the sum of two primes, that
is,

\begin{equation}\label{dois-1}
\forall p \forall q \hspace{8pt} [ p + q \ne m ]
\end{equation}

where $p$ and $q$ belong to the set of prime numbers. The formula
\eqref{dois-1} can be rewritten as

\begin{equation}\label{dois-2}
\forall p \forall q \hspace{8pt} [ (p+1) + q \ne m+1 ]
\end{equation}

The integers $p$ and $q$ are not equal to 2 because if one is
equal to 2 the other must be equal to 2 and this contradicts the
hypothesis made about $m$. The number $p+1$ is an even integer and
it can be expressed as the sum of two odd integers $j$ and $k$,
that is,

\begin{equation}\label{dois-3}
(p+1) + q = j + k + q \ne m+1
\end{equation}

If the ternary GC is true there are three prime numbers $a$, $b$
and $c$ such that

\begin{equation}\label{dois-4}
a + b + c = m+1
\end{equation}

Since that $m+1$ is an odd integer we have the alternatives: ($i$)
$a$, $b$ and $c$ are odd integers or ($ii$) $a$ and $b$ are equal
to 2 and $c$ is an odd integer. Considering the alternative ($i$)
and comparing \eqref{dois-4} with \eqref{dois-2} and
\eqref{dois-3} we obtain a contradiction. Considering the
alternative ($ii$) we have

\begin{equation}\label{dois-5}
3 + c = m
\end{equation}

Comparing \eqref{dois-5} with \eqref{dois-1} we again obtain a
contradiction. The reason of the contradictions is the hypothesis
that ``there is an even integer $m$ that can not be expressed as
the sum of two primes". Therefore if every odd integer $n > 5$ is
the sum of three primes, then every even integer $n > 2$ is the
sum of two primes.

\section{A conditional proof of the ternary GC}
\label{sec-3}

\hspace{22pt} Let $A$ be the well-ordered set of odd integers
greater than 5. Let us denote by $B$ the finite well-ordered
subset of $A$ such that ($i$) each element of $B$ is a value for
which is unknown if the ternary GC is true, $(ii)$ each element of
$A - B$ is a value for which the ternary GC is true and ($iii$)
the set of elements of $A$ less than $\alpha$ is not empty.

\begin{equation}\label{tres-1}
B = \{\alpha, \alpha + 2, \alpha + 4, \ldots , \beta -4, \beta -2,
\beta \}
\end{equation}

If Goldbach's conjecture is true for all even integers less than
$\alpha$ and if any element $n$ of $B$ satisfies an equation of
the form

\begin{equation}\label{tres-2}
n = p + r
\end{equation}

where $p$ is a odd prime number and $r$ belongs to the
well-ordered subset

\begin{equation}\label{tres-3}
\{ 4, 6, 8, \ldots , \alpha - 5, \alpha - 3, \alpha - 1\}
\end{equation}

of $A$, we have that the ternary GC is true for all elements of
$B$.

\hspace{22pt} If the gap between each prime less or equal to
$\beta - (\alpha - 1)$ and its consecutive prime is less than or
equal to $\alpha - 4$, then some $p$ that satisfies \eqref{tres-2}
exists for any element $n$ of $B$. Let us assume that Cram\'er's
conjecture ~\cite{Cramer1} is true, that is,

\begin{equation}\label{tres-4}
\max \limits_{p_n \le k} \ (p_{n+1} - p_n) \sim \log^2 k
\end{equation}

where $k$ is an integer and $p_n$ is the $n$-th element of the
well-ordered set of prime numbers. Substituting $k$ by $\beta$ in
\eqref{tres-4}, we have

\begin{equation}\label{tres-5}
\max \limits_{p_n \le \beta} \ (p_{n+1} - p_n) \sim \log^2 \beta
\end{equation}

Let us consider $\beta$ equal to $10^{7194}$, with this value of
$\beta$ we obtain a maximum gap of

\begin{equation}\label{tres-6}
\sim \Big( \frac{7194}{\log_{10} e} \Big)^2 < 274400000
\end{equation}

\begin{theorem} \label{teorema-2}
If
\begin{equation}\label{tres-7}
\max \limits_{p_n \le \beta} \ (p_{n+1} - p_n) < \log^r \beta
\end{equation}

where
\begin{equation}\label{tres-8}
\log^r \beta = \alpha
\end{equation}

\begin{equation}\label{tres-9}
\alpha > 10 \times \log^2 \beta
\end{equation}

and the ternary GC is true for odd integers less than $\alpha$ and
is also true for odd integers greater than $\beta$, then the
ternary GC is true.

\end{theorem}
\textit{Proof:}

\hspace{22pt} Let us consider $\alpha = 1 \times 10^{16}$
~\cite{Silva} and $\beta = 10^{7194}$ ~\cite{Chen&Wang2}. With
these values we have that the right member of \eqref{tres-9} is
\begin{equation}\label{tres-10}
\approx 2.744 \times 10^9
\end{equation}

and the value of $r$ is
\begin{equation}\label{tres-11}
r = \log_{ ( \log \beta )} \alpha
\end{equation}

\begin{equation}\label{tres-12}
r \approx 3.7921
\end{equation}

Comparing \eqref{tres-10} with $\alpha$ we see that the gap
between each prime less or equal to $\beta - (\alpha - 1)$ and its
consecutive prime is less than or equal to $\alpha - 4$. Therefore
exists some $p$ that satisfies \eqref{tres-2} for any element $n$
of $B$. Considering that Goldbach's conjecture is true for all
even integers less than $\alpha$, we can conclude that if
\eqref{tres-7} holds, then the ternary GC is true.

\hspace{22pt} The statement of theorem ~\ref{teorema-2} assumes in
\eqref{tres-9} that Cram\'er's conjecture is true in the worst
case. With the values used for $\alpha$ and $\beta$, we have a
proof of the ternary GC if the following conjecture is true for
odd integers less than or equal to $\beta$

\begin{equation}\label{tres-13}
\max \limits_{p_n \le \beta} \ (p_{n+1} - p_n) < \log^{3.7921}
\beta
\end{equation}

with $\beta = 10^{7194}$.

\section{ On Cram\'er's model }
\label{sec-4}

\hspace{22pt} In 1943 Selberg ~\cite{Selberg} proved, assuming
Riemann's hypothesis, that

\begin{equation}\label{quatro-1}
\pi ( x + \Phi (x)) - \pi (x) \sim \frac{\Phi (x)}{\log x}
\hspace{22pt} (x \rightarrow \infty )
\end{equation}

is true for almost all $x$ if
\begin{equation}\label{quatro-2}
\frac{\Phi (x)}{\log^{2} x} \rightarrow \infty \hspace{22pt} (x
\rightarrow \infty )
\end{equation}

and, in 1985, Maier ~\cite{Maier} concluded that Selberger's
result is true with exceptions. To try to guess the  number of
primes between $x$  and  $x+y$, Maier  first removed those
integers that have a small prime factor ({\sl following
Eratosthenes}),  and only then did he apply density arguments
({\sl following Gauss}) ~\cite{Granville1}. Maier's result
contradicts what one expected from Cram\'er's model.

\hspace{22pt} Let us consider two integers $x$ and $y$ such that
\begin{equation}\label{quatro-3}
x > y
\end{equation}
and
\begin{equation}\label{quatro-4}
x + y \ggg p
\end{equation}
where $p$ is the greatest prime less than or equal to the square
root of $x+y$. Let us assume that the chance that a given integer
$n$ be a prime is $\frac {1}{\log n}$. Let $q$ be a prime less
than or equal to $p$ and let $F(n,q)$ be a function such that
\begin{equation}\label{quatro-5}
\frac {F(n,q)}{\log n}
\end{equation}
is the chance that a given integer $n$ belonging to $(x,x+y]$ be
prime after crossing out of the interval those integers that are
divisible by the primes less than or equal to $q$. For $q=p$ we
have
\begin{equation}\label{quatro-6}
\frac {F(n,p)}{\log n} = 1
\end{equation}
or
\begin{equation}\label{quatro-7}
F(n,p) = \log n
\end{equation}

In accordance with Friedlander and Granville
~\cite{Friedlander&Granville1} and Maier ~\cite{Maier} the
probability of a `randomly chosen' integer $n$ belonging to
$(x,x+y]$ being prime, given that it has no prime factors $\le z$
and being $z$ a `small' prime, is
\begin{equation}\label{quatro-8}
\frac {1}{ \Big\{ \prod_{s \le z} \Big( 1 - \frac {1}{s} \Big)
\Big\} \log n}
\end{equation}

Considering \eqref{quatro-4} let $z$ be equal to $p$. Since that
\begin{equation}\label{quatro-9}
\frac {1}{ \Big\{ \prod_{s \le p} \Big( 1 - \frac {1}{p} \Big)
\Big\} } \ne \log n
\end{equation}
we have a contradiction between \eqref{quatro-7}  and
\eqref{quatro-9}. The probability \eqref{quatro-8} is crucial in
Maier's work ~\cite{Maier}.

\end{document}